\documentclass[10pt,conference]{IEEEtran}
\IEEEoverridecommandlockouts
\usepackage{amssymb,amsmath}
\usepackage{graphicx}
\usepackage{epsfig}
\usepackage{grffile}
\usepackage{balance}
\pdfoutput=1
\begin{document}

\title{Fair Allocation of Backhaul Resources in Multi-Cell MIMO Co-operative Networks }
\author{Hazem M. Soliman, Omar A. Nasr and Mohamed M. Khairy \\
Center for Wireless Studies, Faculty of Engineering\\ 
Cairo University\\
\thanks{This work was funded by the Egyptian National Telecommunications Regulatory Authority (NTRA).}
}
\date{}
\maketitle

\thispagestyle{empty}

\begin{abstract}
In this paper the problem of allocating the limited backhaul bandwidth among users in Multi-cell MIMO cooperative
networks is considered. We approach the problem from both the sum-rate and fairness perspectives. First, we show that
there are many allocations that can provide near maximum sumrate while varying significantly in fairness, which is assessed
through the mean versus variance criteria. Second, Two novel schemes that achieve near maximum sum-rate while offering
fair allocation of the backhaul bandwidth among users are proposed: the Equal Signal-to-Interference ratio (SIR) and the
Equal Interference schemes. Simulation results show that, for the same mean rate among users, the proposed schemes can achieve
more fairness when compared to the conventional scheme, which gives all users the same share of bandwidth. Moreover, we show
that the Equal SIR scheme can achieve zero variance among users in a wide range of backhaul bandwidths while keeping
very close to maximum sum rate. This is the most fair solution that can be used in Multi-cell MIMO in that range of backhaul
bandwidths.

\end{abstract}


\section{Introduction}
\label{sec:intro}

Interference is one of the key challenges facing  future wireless communication systems. The conventional technique to deal with interference is to limit the re-usability of resources (time, frequency, code,...) to introduce some kind of orthogonality between users. The more recent approach tries to make use of interference or at least coordinates it through the use of cooperative transmission \cite{Gesbert2010}. Multi-Cell MIMO or Network MIMO is a new technology for cellular base-stations that mitigates interference by coordinating base-stations transmission. 3GPP LTE-A and IEEE 802.16m have recently chosen Network MIMO as a means to  increase the cell-edge and system throughput \cite{3GPPOctober2008} \cite{HuangJanuary2008}. 

A fundamental challenge in Multi-Cell MIMO networks is the issue of limited backhaul bandwidth. Multi-Cell processing or full co-operation among base-stations requires the exchange of full channel state information (CSI) and data transmitted to users among all base-stations, this requires a very high-speed backhaul. Another type of collaboration called Interference Coordination \cite{Gesbert2010} requires the exchange of CSI only to perform  some form of coordinated beamforming. Several attempts to reduce backhaul requirements through distributed cooperation, statistical CSI exchange or clustered cooperation have been proposed \cite{ZakhourSpring2009} \cite{Zakhour2010} \cite{ZhangAPRIL2009}.

Another approach to reduce backhaul load is to do co-operation only for a selected subset of users according to a criteria that selects only the deserving users \cite{Marsch2008} \cite{Marsch2009}  \cite{Lei2009}. Authors of those papers have mainly considered data sharing, focusing on sum-rate   and not fairness. Up to our knowledge,  the problem of analytical optimization of sharing  the limited backhaul between users for coordinated beamforming in Network MIMO has not been considered before.

The installation of Network MIMO based base-stations introduces  new resources in the network. One type of these resources is the backhaul bandwidth.In any resource allocation problem, there is always a trade-off between global performance, best represented by sum-rate, and individual performance, best represented by fairness. Although  fairness is usually studied in higher layers, the emerging cross-layer design concepts have encouraged the study of fairness at physical layer. A comprehensive study of fairness in wireless communications both in physical and MAC layers, was performed in \cite{Calvo2004}.  The used criteria is to measure the mean vs. variance among the users. It was shown that most of the time whenever global performance, mean, is maximized, the difference between maximum and minimum, variance, increases, meaning worse fairness. 

 In this paper we propose two schemes to distribute the backhaul bandwidth  ensuring fairness among users, measured by the mean vs. variance metric. We first derive an analytical model for the interference in coordinated beamforming under a limited backhaul bandwidth assuming only quantized CSI is shared but no data sharing. We then study the problem of allocating unequal shares of  the backhaul bandwidth to users from two different perspectives. First we study the sum-rate problem, where we analytically show that  the optimization of sum-rate provides insignificant performance improvement compared to the conventional scheme which equally allocates backhaul bandwidth among users. Secondly, and based on the previous important observation, we study the problem from the fairness perspective and propose two schemes to distribute the backhaul bandwidth, based on  Equal Signal-to-Interference ratio (SIR) and  Equal Interference criteria. We then show, through analysis and simulation, that these two approaches are  able to provide the same  sum-rate while ensuring fairness among users compared to the conventional scheme. 
 
The paper is organized as follows, Section \ref{sec:sysmodel} introduces the system  model. Section \ref{sec:sumrate} describes the problem formulation, solution and fairness discussion. Section \ref{sec:simres} shows the simulation results.
Finally section \ref{sec:conclusion} concludes the paper.

\section{System  Model}
\label{sec:sysmodel}

We consider a Wyner type \cite{WynerNov.1994}, two base stations, $N$-user per cell MIMO downlink system. Each base station has $M$ antennas, while users are each equipped with a single antenna. Channel is taken from the Zero-mean Circularly-symmetric Complex-Gaussian model $(ZMCSCG)$\cite{Paulraj2003}. We also take the simple yet efficient Zero-Forcing precoder, its columns are  normalized to obey power limit condition. Hence the received signal may be expressed as
\begin{equation}
\mathbf{y} = \mathbf{H}\mathbf{W}\mathbf{x} + \mathbf{n}
\end{equation}
where $\mathbf{x}$ is the transmitted signal, $\mathbf{y}$ is the received signal vector, $\mathbf{H}$ is the channel matrix, $\mathbf{W}$ is the linear beamformer used, and $\mathbf{n}$ is the noise vector.
The instantaneous received signal to interference and noise ratio $\mathit{(SINR^t_i)}$ for a general beamformer for user $i$ in the first cell is

\begin{equation}
\displaystyle \mathit{SINR}^t_i = \frac{ P_{1,i} \left| \mathbf{h}_{1,i} \mathbf{w}_{1,i} \right| ^2}{ \sigma^2 +\sum\limits_{\substack{ 1<j<N \\ i \neq j}} P_{1,j} \left| \mathbf{h}_{1,i}^H \mathbf{w}_{1,j} \right| ^2 + \sum\limits_{\substack{ 1<j<N }} P_{2,j} \left| \mathbf{h}_{2,i}^H \mathbf{w}_{2,j} \right| ^2 }
\label{eq:sinr}
\end{equation}
where $P_{k,i}$ is the power transmitted from base station $k$ intended to user $i$ and $\sigma^2$ is the noise power.
Assuming each user has a rate $r_i$ (in bits/symbol/Hz), the sum-rate for $N$-users in each cell is given by
\begin{equation}
\sum_{1<i<N} r_i =\sum_{1<i<N} \log_2{ \left( 1 + \mathit{SINR}^t_i \right) }
\end{equation}

Each base station is assumed to have perfect knowledge of its own users' channels. These channels should then be conveyed through the backhaul to the other base station to do coordinated beamforming. However, due to limited backhaul, only quantized versions of the channels may be exchanged between the coordinating base-stations. Uniform quantization is assumed. Each base station uses the perfectly known channels of its own users, and the quantized channels of the other cell users to design a beamforming matrix $\mathbf{W}$. Clearly, even with coordinated Zero-forcing beamforming, the users whose channels were quantized will suffer from interference due to quantization. This quantization interference is what is left from the multi-cell interference and is still  the major limit for the system performance. The problem of how to allocate the quantization bits among users is considered in the following Section.


\section {Sum-Rate vs. Fairness Analysis}
\label{sec:sumrate}

\begin{figure}[h]
	\centering
		\includegraphics[width=0.45\textwidth, height=0.45\textwidth]{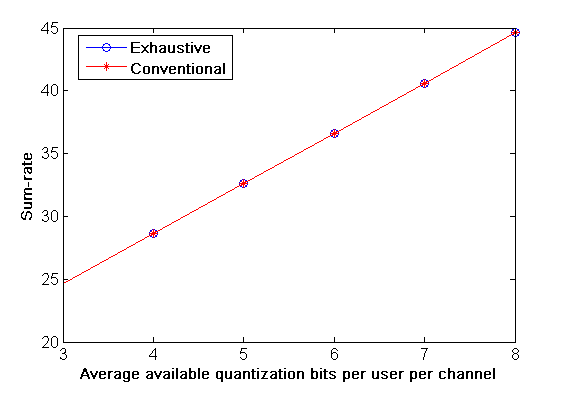}
	\caption{Sum-Rate achieved using exhaustive search and conventional scheme $\{P_{1,i}=P_{2,i}=10 \;\; \forall i=1:N, N=4, M=8\}$ }
	\label{fig:B_W0}
\end{figure}

In Fig. \ref{fig:B_W0}  the sum-rate using two schemes is shown. The first scheme is the conventional scheme which assigns the same number of backhaul bits for all users, and the second  is an exhaustive search, which optimally allocates quantization bits among users to maximize sum-rate. Surprisingly,  the optimum scheme provides insignificant  performance improvement over the conventional scheme. Although several schemes might have almost the same sum-rate, they can have entirely different fairness performance as will be shown in Section \ref{sec:simres}.
In the following subsections, we explain this phenomena, then we turn our attention to fairness where we propose two algorithms that provide more fairness while keeping a  close-to-optimum sum-rate.

\subsection{Problem Formulation}
Maximization of  the sum rate under a constrained backhaul, which is used to convey quantized versions of the users channels is considered in this Section. First, we derive an analytical formulation for the average received $SINR_i$ using a beamformer based on the quantized-channel, then we formulate and solve the sum-rate maximization problem. 
It should be noted that the interference from the same cell users equals zero due to using zero-forcing beamforming based on perfectly known channels. The received multi-cell interference $\mathit{I}_{i}$ at user $i$ in the first cell  is
\begin{equation}
\mathit{I}_{i} = \sum_{\substack{ 1<j<N }} \mathit{P}_{2,j} \left| \mathbf{h}_{2,i}^H \mathbf{w}_{2,j} \right| ^2 
\end{equation}
the unquantized channel vector $\mathbf{h}_{2,i}$ may be writtem in terms of the quantized channel vector $\mathbf{hq}_{2,i}$ using the simple formula
\begin{equation}
\mathbf{h}_{2,i}=\mathbf{hq}_{2,i} + \mathbf{nq}_{i}
\end{equation}
where  $\mathbf{nq}_i$ is the quantization noise vector corresponding to user $i$ channels. So
\begin{equation*}
\mathit{I}_{i} = \sum_{\substack{ 1<j<N }} \mathit{P}_{2,j} \left| (\mathbf{hq}_{2,i} + \mathbf{nq}_{i})^H \mathbf{w}_{2,j}\right| ^2 
\end{equation*}
\begin{align*}
\mathit{I}_{i} =& \sum_{\substack{ 1<j<N }} \mathit{P}_{2,j}\left[ \left| \mathbf{hq}_{2,i}^H \mathbf{w}_{2,j} \right| ^2 + \left| \mathbf{nq}_{i}^H \mathbf{w}_{2,j} \right| ^2 \right. \\
& \left. + \left(\mathbf{hq}_{2,i}^H \mathbf{w}_{2,j} \right)^*\left(\mathbf{nq}_{i}^H \mathbf{w}_{2,j}\right) + \left(\mathbf{hq}_{2,i}^H \mathbf{w}_{2,j} \right)\left(\mathbf{nq}_{i}^H \mathbf{w}_{2,j}\right)^* \right]
\end{align*}
the first term $ \left| \mathbf{hq}_{2,i}^H \mathbf{w}_{2,j} \right| ^2$ is independent on the quantization noise and depends only on the beamformer type. It also equals zero for the Zero-forced beamformer we use throughout the  paper.
So the interference due to quantization only, $\mathit{I}_{qi}$,  equals the total interference $\mathit{I}_{i}$ and is given by
\begin{equation}
\begin{split}
\mathit{I}_{qi} =& \sum_{\substack{ 1<j<N }} \mathit{P}_{2,j}\left[ \left| \mathbf{nq}_{i}^H \mathbf{w}_{2,j} \right| ^2 + \left(\mathbf{hq}_{2,i}^H \mathbf{w}_{2,j} \right)^*\left(\mathbf{nq}_{i}^H \mathbf{w}_{2,j}\right) \right. \\
& \left. + \left(\mathbf{hq}_{2,i}^H \mathbf{w}_{2,j} \right)\left(\mathbf{nq}_{i}^H \mathbf{w}_{2,j}\right)^* \right]
\end{split}
\end{equation}
taking mean over channel realizations
\begin{equation}
\mathbf{E}\left \{ \mathit{I}_{qi} \right \}=\mathbf{E}\left \{ \sum_{\substack{ 1<j<N }} \mathit{P}_{2,j}\left( \left| \mathbf{nq}_{i}^H \mathbf{w}_{2,j} \right| ^2 \right) \right \}
\end{equation}
 as the  terms 
 \begin{equation}
 \mathbf{E} \{\left(\mathbf{hq}_{2,i}^H \mathbf{w}_{2,j} \right)^*\left(\mathbf{nq}_{i}^H \mathbf{w}_{2,j}\right) + \left(\mathbf{hq}_{2,i}^H \mathbf{w}_{2,j} \right)\left(\mathbf{nq}_{i}^H \mathbf{w}_{2,j}\right)^* \}
 \end{equation}
  will equal zero due to the zero-mean of the quantization noise. Hence
\begin{equation*}
\mathbf{E}\{ \mathit{I}_{qi} \}= \sum_{\substack{ 1<j<N }} \mathit{P}_{2j} \mathbf{E} \left \{\left| \sum_{\substack{ 1<k<M }} {nq}_{i,k} {w}_{2,j,k} \right| ^2 \right \}
\end{equation*}
where ${nq}_{i,k}$,${w}_{2,j,k}$ are the $k^{th}$ elements of $\mathbf{nq}_{i}$ and $\mathbf{w}_{2,j}$ respectively.
\begin{equation*}
\mathbf{E}\{ \mathit{I}_{qi} \}= \sum_{\substack{ 1<j<N }} \mathit{P}_{2,j}  \sum_{\substack{ 1<k<M }}  \mathbf{E}\{ \left| {nq}_{i,k} {w}_{2,j,k} \right| ^2 \}
\end{equation*}
\begin{equation*}
\mathbf{E}\{ \mathit{I}_{qi} \}= \sum_{\substack{ 1<j<N }} \mathit{P}_{2,j}  \sum_{\substack{ 1<k<M }}  \mathbf{E}\{\left| {nq}_{i,k}\right|^2 \}   \mathbf{E}\{\left|{w}_{2,j,k} \right| ^2 \}
\end{equation*}
\begin{equation*}
\mathbf{E}\{ \mathit{I}_{qi} \}= \sum_{\substack{ 1<j<N }} \mathit{P}_{2,j} \mathit{Q}_i  \sum_{\substack{ 1<k<M }} \left| {w}_{2,j,k} \right| ^2 
\end{equation*}
and as we assumed before that beamformer vectors are normalized $\left( \sum_{\substack{ 1<k<M }} \left| {w}_{2,j,k} \right|^2= 1 \right)$, then
\begin{equation}
\mathbf{E}\{ \mathit{I}_{qi} \} = \sum_{ 1<j<N }{ \mathit{P}_{2,j} \times \mathit{Q}_i}
\label{eq:intf}
\end{equation}
where $Q_i$ is the quantization noise and is given by \cite{Haykin2001}
\begin{equation}
\mathit{Q}_i = \mathbf{E}\{ \left| {nq}_{i,k}\right|^2 \}= \frac{\text{const}}{2^{-2 \times l_i}}
\label{eqn:qeqn}
\end{equation}
where $l_i$ is the number of quantization bits per channel given to user $i$.
\begin{figure}[h]
	\centering
		\includegraphics[width=0.45\textwidth, height=0.45\textwidth]{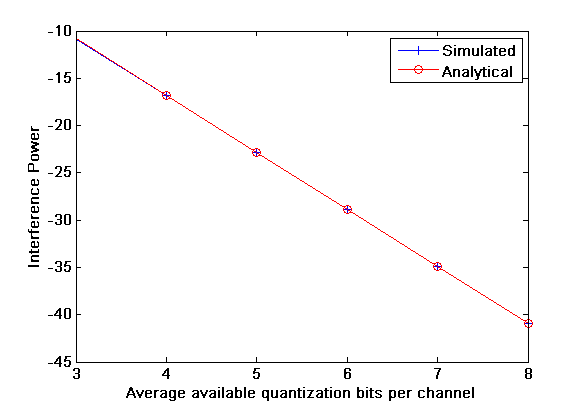}
	\caption{Analytical \& Simulated interference  $\{P_{1,i}=P_{2,i}=10 \;\; \forall i=1 : N, N=8\}$ }
	\label{fig:B_W1}
\end{figure}
The simulated interference and (\ref{eq:intf}) are shown in Fig.  \ref{fig:B_W1} where it is clear that the two curves are matching. From (\ref{eq:sinr}) and (\ref{eq:intf}), the average $\mathit{SINR}_i$ may be written as
\begin{equation}
\begin{split}
\mathit{SINR}_i = \mathbf{E} \{\mathit{SINR}^t_i\} &=  \frac{ P_{1,i} \mathbf{E} \{  \left| \mathbf{h}_{1,i} \mathbf{w}_{1,i} \right| ^2 \}}{ \sigma^2  + \sum\limits_{ 1<j<N }{ \mathit{P}_{2,j} \times \mathit{Q}_i}}\\
&= \frac{ P_{1,i} }{ \sigma^2  + \sum\limits_{ 1<j<N }{ \mathit{P}_{2,j} \times \mathit{Q}_i}}
\end{split}
\end{equation}
as $\mathbf{E} \{  \left| \mathbf{h}_{1,i} \mathbf{w}_{1,i} \right| ^2 \}=M-N+1$ from \cite{Calvo2004} for the zero forcing beamformer. This value is the same for all users and is omitted for convenience. For the case of uniform power allocation in the second base-station
\begin{equation}
\mathit{SINR}_i = \frac{ P_{1,i} }{ \sigma^2  + \mathit{P}_{2,i} \times \mathit{N} \times \mathit{Q}_i }
\label{eqn:sinreqn2}
\end{equation}

Note that $\mathit{P}_{2,i}$ is  the value of the power  received at user $i$ in the first cell when transmitted from the second base-station intended to its users. Although the average $\mathit{SINR_i}$ is  independent on the channel, the values of $P_{1,i}$ and $P_{2,i}$ vary greatly according to the user position due to the path loss. The main goal of this paper is to allocate bits among users according to their position rather than their instantaneous fading channel. The path-loss model will be introduced at the end of this section.


\subsection{Sum-Rate}
The sum-rate maximization problem is  formulated  as
\begin{equation}
\begin{split}
\text{max} &    \;\;\;\;\;  \sum\limits_{\substack{1<i<N}} r_i\\
\text{s.t.} &   \;\;\;\;\; \sum_{\substack{1<i<N}} l_i = \mathit{D}
\end{split}
\end{equation}
The Lagrange multiplier formulation of the above problem is
\begin{equation}
\mathbf{J}= \sum_{\substack{1<i<N}} \log_2{ \left(1 + \mathit{SINR}_i \right)} + \mu  \left( \mathit{D} - \sum_{\substack{1<i<N}} l_i \right)
\end{equation}
where $D$ is the total number of backhaul bits and $\mu$ is the Lagrange multiplier, differentiating and equating with zero yields the following condition
\begin{equation}
\frac{ \mathit{I}_{qi} \times \mathit{P}_{1,i} }{ \left( \mathit{I}_{qi} + \sigma^2 \right) \times \left( \mathit{P}_{1,i} + \mathit{I}_{qi} + \sigma^2 \right)} = \mathbf{A}
\label{eq:cond}
\end{equation}
where $\mathbf{A}$ is some constant. Next, we divide the system operation into 2 regions as follows:\\
\subsubsection{Region 1}
\begin{equation}
\sigma^2  \ll \mathit{I}_{qi} \ll \mathit{P}_{1,i},      \;\;\;\;\;         \forall i=1 : N
\end{equation}
In region 1, a good approximation for (\ref{eq:cond}) may be as follows  
\begin{equation}
\frac{\mathit{I}_{qi} \times \mathit{P}_i}{\mathit{I}_{qi} \times \mathit{P}_i} = \mathsf{1} = \mathbf{A''}
\end{equation}
 this result clearly shows that the condition for sum-rate maximization is trivial and can be easily achieved with any scheme that gives reasonable number of bits to the edge users which represents the lowest $\mathit{SINR}$ or the largest interference. Examples for these schemes are the conventional scheme or the more fair schemes proposed  later in this section.

\subsubsection{Region 2}
\begin{equation}
\sigma^2  \approx \mathit{I}_{qi} \ll \mathit{P}_{1,i},      \;\;\;\;\;         \forall i=1 : N
\end{equation}
A good approximation here is 
\begin{equation}
\frac{\mathit{I}_{qi} \times \mathit{P}_i}{ \left( \mathit{I}_{qi} + \sigma^2 \right) \times \mathit{P}_i} = \mathbf{A}
\end{equation}
or
\begin{equation}
\mathit{I}_{qi} = \mathbf{A'''}
\label{eqn:reg3eqn}
\end{equation}
which means equal interference for all users. However in region 2, similar to the famous power water-filling problem where the water-filling and the uniform power allocation are the same for high $SNR$ \cite{Haykin2001}, we argue that  in region 2 the number of bits is  large enough to make all schemes, the conventional scheme and the two water-filling like schemes proposed next, approach each other as will be shown from the simulation results.

\subsection{Fairness}
Based on our discussion of Fig. \ref{fig:B_W0}, we argue that  the sum-rate maximization condition is  a very relaxed one and can be easily achieved using  several schemes. That is why we choose to turn our attention  to fairness issues. As  mentioned in the introduction and in \cite{Calvo2004}, global performance is usually penalized whenever we want better fairness. However, our important contribution in this paper is that in some cases, where the conditions for global performance is relaxed as we showed before, clever allocation of the backhaul bandwidth can be utilized to provide a much better fair performance while keeping the global performance, the sum-rate, almost the same. Fair schemes can by obtained  by solving a max-min optimization problem for a specific system performance criteria. The optimum solution for a max-min problem is  water-filling  \cite{RadunovicOctober2007}.  The two proposed fair schemes are obtained by using the SIR and interference as subject of max-min optimization. These schemes are    ({\it{Equal SIR}}) and {\it{(Equal Interference)}}. These two fair schemes provide better fair conditions while preserving the global performance as well.
\subsubsection{ Equal Signal-to-Interference-ratio (SIR)}
in order to provide equal SIR for all users, we solve the following problem
\begin{equation}
\begin{split}
\text{max min}&   \;\;\;\;\;  \frac{P_{1,i}}{I_{qi}}\\
\text{s.t.}&   \;\;\;\;\; \sum_{\substack{1<i<N}} l_i = \mathit{D}
\end{split}
\label{eq:eqsirprob}
\end{equation}
from (\ref{eqn:qeqn}) and (\ref{eqn:sinreqn2}), the solution for (\ref{eq:eqsirprob}) is a water-filling like equation for allocating bits among users as follows\\
\begin{equation}
l_i= \mathit{a} + 0.5 \times \log_2{\left( \frac{ \mathit{P}_{2,i}}{ \mathit{P}_{1,i}}\right)}
\end{equation}
where $a$ is a constant depending on the total number of bits and can be found by solving the equation
\begin{equation}
\sum_{1<i<N} l_i = \sum_{1<i<N} \left ( \mathit{a} + 0.5 \times \log_2{\left( \frac{ \mathit{P}_{2,i}}{ \mathit{P}_{1,i}}\right)} \right ) = \mathit{D}
\end{equation}
where $D$ is the total number of bits available in the backhaul.

\subsubsection{ Equal Interference}
in order to provide equal interference among users, we solve the following problem
\begin{equation}
\begin{split}
\text{max min}&    \;\;\;\;\;  I_{qi}\\
\text{s.t.}&   \;\;\;\;\; \sum_{\substack{1<i<N}} l_i = \mathit{D}
\end{split}
\label{eq:eqiprob}
\end{equation}
from (\ref{eqn:qeqn}), (\ref{eqn:sinreqn2}) and (\ref{eqn:reg3eqn}), the solution for (\ref{eq:eqiprob}) is a water-filling like equation for allocating bits among users  as follows
\begin{equation}
l_i= \mathit{a} + 0.5 \times \log_2{\left( { \mathit{P}_{2,i}}\right)}
\end{equation}
where $a$ is a constant depending on the total number of bits and can be found by solving the equation
\begin{equation}
\sum_{1<i<N} l_i = \sum_{1<i<N} \left ( \mathit{a} + 0.5 \times \log_2{\left( { \mathit{P}_{2,i}}\right)} \right ) = \mathit{D}
\end{equation}
where $D$ is the total number of bits available in the backhaul.\\
We consider the path-loss model of \cite{KARAKAYALI2006}, in which the power received at distance $d$ can be expressed as
\begin{equation}
P(d) = P_o \times \mathit{k} \times \left(\frac{d}{d_o}\right)^{-\gamma}
\end{equation}
 where  $ d_{o}$  is the reference distance, $k$ is the loss at this reference distance and $\gamma$ is the path-loss exponent. Using this model, $\mathit{SINR}_i$ is expressed as
\begin{equation}
\mathit{SINR}_i = \frac{\mathit{P} \times \mathit{k} \times \left( \frac{d_{1,i}}{d_o} \right)^{-\gamma}}{\sigma^2 + \mathit{P} \times \mathit{k} \times \mathit{N} \times \left(\frac{d_{2,i}}{d_o} \right)^{-\gamma} \times \mathit{Q}_i}
\end{equation}
Also expressions for the number bits per user in the Equal SIR scheme is given by 
\begin{equation}
l_i= \mathit{a} + 0.5 \times \log_2{\left( \frac{ \mathit{d}_{2,i}}{ \mathit{d}_{1,i}}\right)^{-\gamma}}
\end{equation}
and in the Equal Interference scheme is given by
\begin{equation}
l_i= \mathit{a} + 0.5 \times \log_2{\left( { \mathit{d}_{2,i}}\right)^{-\gamma}}
\end{equation}\\

\section{Simulation Results}
\label{sec:simres}

\begin{figure}[!]
	\centering
		\includegraphics[width=0.5\textwidth, height=0.5\textwidth]{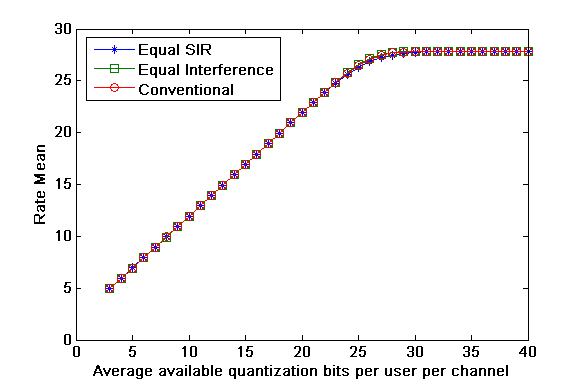}
	\caption{Rate Mean versus number of quantization bits}
	\label{fig:C_W1}
\end{figure}

\begin{figure}[!]
	\centering
		\includegraphics[width=0.5\textwidth, height=0.5\textwidth]{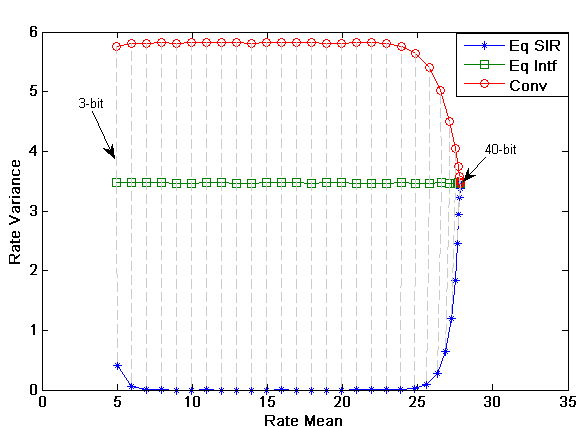}
	\caption{Rate variance vs. Rate mean with vertical lines represents a constant number of quantization bits}
	\label{fig:C_W3}
\end{figure}

The simulation is done for a two-base stations, 8 users per cell. Backhaul bandwidth is increased from 3 to 40 bits on average per user. The reference distance $d_o$ is assumed to be 1600m and is also considered the cell radius. Path loss exponent is 3.8. Users per cell are allocated every 200m starting from 200m to 1600m. Base station power is equal to 10 watt. Two important parameters are used to measure system performance: rate mean, which is defined as the average rate over all users’ rates, and rate variance, which is the variance among different rates of the users. As mentioned earlier, less variance for the same mean indicates better fairness. Fig. \ref{fig:C_W1} plots the rate mean vs. average number of available bits per user per channel. It is clear from the figure that all three schemes achieve the same mean, the explanation for this result was given in Section \ref{sec:sumrate}. Fig. \ref{fig:C_W3}  is the combined plot of the mean vs. variance which shows that our proposed schemes can achieve a much less variance while preserving the same mean. Also in  Fig. \ref{fig:C_W3} we can clearly see the two regions of operation of the system. The first region representing moderate interference  ranges from $3 \sim 25$ bits and this is where our schemes perform at their best from the fairness point of view. The second region  ranges  from 26 bits onward and this is where the three schemes converge. The second region ends with all schemes approaching each other and this can be considered as  the infinite backhaul point. Finally we may have three notes. First is that although the first scheme, Equal SIR, is much better than the second one, Equal Interference, the second algorithm is proposed because it is  the one that all other algorithms converge to at region 2. Second note is that the Equal SIR scheme can achieve the lowest possible value for variance which is zero. Last note is that asymptotically, when we approach the infinite backhaul point, no one can achieve zero rate variance.  In this region, all users suffer from no interference because channels are transmitted through the backhaul with no quantization. Consequently, because different users are located at different distances from the basestations, they will receive different powers. Hence, the rates of different users will not be the same, and the rate variance will no longer equals zero. 

\section{conclusion}
\label{sec:conclusion}

In this paper we studied the problem of how to best allocate the backhaul bandwidth among users in a multi-cell MIMO system. We solved the problem of exchanging CSI between two base stations for coordinated beamforming. The first approach was to maximize sum-rate where we showed that the conditions for its optimization are very relaxed. We then turned our attention towards fairness where we proposed two schemes, the Equal SIR and Equal Intereference. We then showed through simulations how the proposed schemes can achieve a much less variance while keeping the mean rate performance as is.

\balance

\bibliographystyle{IEEEtrannames}
\bibliography{hazem2010}

\end{document}